\documentclass[letterpaper, 10 pt, conference]{ieeeconf}  

\IEEEoverridecommandlockouts                              
\usepackage{times}
\usepackage{subfig}
\usepackage{graphicx}
\usepackage{amsmath}
\usepackage{amsthm}
\usepackage{amssymb}
\usepackage{amsfonts}
\usepackage{mathtools}
\usepackage{xcolor}
\usepackage{booktabs}
\usepackage{array}
\usepackage{cite}
\usepackage{adjustbox}
\let\labelindent\relax
\usepackage{enumitem}
\usepackage[font=footnotesize,labelfont=bf]{caption}
\usepackage{tikz-network}

\DeclareMathOperator{\tr}{tr}

\newcommand{\PLI}{P\L{}I\;}
\newcommand{\gPLI}{\textit{g}P\L{}I\;}
\newcommand{\satPLI}{\textit{sat}P\L{}I\;}

\newtheorem{theorem}{Theorem}
\newtheorem{lemma}{Lemma}
\newtheorem{proposition}{Proposition}
\newtheorem{corollary}{Corollary}

\newtheorem{definition}{Definition}

\newtheorem{remark}{Remark}

\newtheorem{problem}{Problem}

\usepackage{tikz-network}
\tikzset{
    block/.style = {draw, fill=white, rectangle, minimum height=.75cm, minimum width=1.25cm},
    Circle/.style = {draw, thick, fill=white, circle, minimum height=.75cm, minimum width=1.25cm},
    tmp/.style  = {coordinate}, 
    sum/.style= {draw=white, fill=white, circle, minimum height=.75cm, minimum width=1.25cm},
    input/.style = {coordinate},
    output/.style= {coordinate},
    pinstyle/.style = {pin edge={to-,thin,black}
    }
}

\title{\LARGE \bf
On the (almost) Global Exponential Convergence of the Overparameterized Policy Optimization for the LQR Problem
}

\author{Moh Kamalul Wafi$^1$, Arthur C. B. de Oliveira$^1$, and Eduardo D. Sontag$^{1,2}$
\thanks{$^1$Department of Electrical and Computer Engineering, Northeastern University, USA.
{\tt\footnotesize	\{wafi.m, a.castello\}@northeastern.edu}.}
\thanks{$^2$Department of BioEngineering and affiliated with the departments of Chemical Engineering and Mathematics, Northeastern University, USA. 
{\tt\footnotesize	e.sontag@northeastern.edu}.}
}
\begin{document}

\maketitle
\thispagestyle{empty}
\pagestyle{empty}

\begin{abstract}
In this work we study the convergence of gradient methods for nonconvex optimization problems -- specifically the effect of the problem formulation to the convergence behavior of the solution of a gradient flow. We show through a simple example that, surprisingly, the gradient flow solution can be exponentially or asymptotically convergent, depending on how the problem is formulated. We then deepen the analysis and show that a policy optimization strategy for the continuous-time linear quadratic regulator (LQR) (which is known to present only asymptotic convergence globally) presents almost global exponential convergence if the problem is overparameterized through a linear feed-forward neural network (LFFNN). We prove this qualitative improvement always happens for a simplified version of the LQR problem and derive explicit convergence rates for the gradient flow. Finally, we show that both the qualitative improvement and the quantitative rate gains persist in the general LQR through numerical simulations.
\end{abstract}
\allowdisplaybreaks

\section{Introduction}\label{sec:intro}

Motivated by the rise in popularity of machine learning methods and neural networks, research on gradient methods and policy optimization and their application have been the focus of much attention in recent years \cite{menon2024geometry,pasande2024stochastic,yuksel2025robust,ahmad2024accelerating,hu2023toward,zhang2023revisiting}. The intuitive idea of finding the minimum of a function by ``traveling along'' its direction of steepest descent stood the test of time, with many recent studies \cite{arora2019implicit,chitour2018geometric,cui2024small,de2024remarks,ArthurAutomatica25} exploring inherent properties that justify the strength of gradient methods.

Despite this success, gradient methods are not guaranteed to solve (or even converge for) every problem, and additional assumptions are usually required. In \cite{karimi2016linear} the authors provide an overview of different assumptions in optimization theory, as well as their relationship. Of particular note is the assumption of convexity, which results in its own field of research in optimization theory \cite{boyd2004convex}.

Alternatively, the Polyak-\L{}ojasiewicz inequality (\PLI) \cite{polyak1963gradient,karimi2016linear} provides a more general guarantee of convergence for solutions of gradient methods, being one of the most important assumptions in nonconvex optimization. Despite that, however, some problems can be shown to not satisfy any \PLI globally, but still be convergent and optimal under gradient methods.

An interesting concrete example of this can be found when studying policy optimization formulations for the linear quadratic regulator (LQR) problem \cite{cui2024small,sontag_remarks_2022,fatkhullin2021optimizing,watanabe2025revisiting,fazel2018global,sun2021learning,hu2023toward,mohammadi2021lack,mohammadi_convergence_2022}. In \cite{hu2023toward} the authors showed that a policy optimization method for the discrete-time LQR problem can be shown to always satisfy a \PLI globally, however existing results for the continuous-time version of the problem \cite{fatkhullin2021optimizing,cui2024small} were only able to show that a \PLI holds on any level-set of the cost function. 

At first glance, the distinction might appear irrelevant, since properness of the LQR cost function imply that any initialization is contained in some of its levelset. However, this distinction motivated follow-up research \cite{cui2024small,ArthurCDC25}. In \cite{cui2024small} the authors first showed its importance by proving distinct robustness properties for problems satisfying either condition. In the same paper the authors also characterized a stronger ``type'' of \PLI inequality for the continuous-time LQR policy optimization problem. In \cite{ArthurCDC25} we formally introduced different types of \PLI and illustrated their relationship and different convergence guarantees. Of particular interest to this paper is the relationship between the usual \PLI (named global \PLI or \gPLI in \cite{ArthurCDC25}) and the saturated \PLI (\satPLI). In \cite{ArthurCDC25} we show that any cost satisfying a \gPLI satisfies a \satPLI, however the converse is not true. We further prove that gradient flows of problems satisfying only a \satPLI present qualitatively worse convergence properties than those satisfying a \gPLI. This illustrates the importance in the distinction between these concepts, and why one might be interested in their relationship.

In this paper, we focus on the \satPLI and \gPLI and their dependence on the problem formulation. Formulation-based acceleration of gradient methods is a well-known phenomenon in the literature, with linear feedforward neural networks (LFFNN) having results in the literature on initialization-dependent quantitative improvements on the convergence rate of gradient methods \cite{ArthurCDC23,min2021explicit,min2023convergence}. However the question we investigate in this paper regards a possible qualitative improvement on the convergence of solutions.

We begin by formally introducing the concepts in section \ref{sec:TheoreticalFramework}, and then illustrate that a simple reparameterization of the form $g(x):=f(x^2)$ can convert a function $f(x)$ that satisfies only a \satPLI (and thus presents only asymptotic convergence to gradient methods) into a function $g(x)$ that can provably be shown to satisfy a \gPLI almost everywhere. In section \ref{sec:LQRPolicyOptimization}, we return to the policy optimization continuous-time LQR as a testbed, and explore the effects of using a LFFNN to the convergence of solutions of a gradient flow. We prove that, for a simplified case of the problem, the overparameterized problem formulation (i.e. using the LFFNNs) improves the convergence of solutions from linear–exponential to exponential, with the trade-off of adding a saddle-point to the dynamics, which can trap or slow-down some solutions. We then provide, in section \ref{sec:NumericalResults}, numerical evidence of this same observation holding for general LQR problem, illustrating the properties discussed theoretically in the previous section. We finish the paper in section \ref{sec:Conclusion-FutureWork} with an overview of our results and some possible future directions of work related to this paper.

\section{Theoretical Framework}\label{sec:TheoreticalFramework}

\subsection{Preliminaries}
Let $\mathcal{X}\subseteq\mathbb{R}^n$ be open and $\mathcal{L}:\mathcal{X}\to\mathbb{R}$ be real-analytic, bounded below, and proper (i.e. coercive or radially unbounded) loss/cost function. We consider the optimization problem
\begin{equation}
    \min\nolimits_{k \in \mathcal{X}} \; \mathcal{L}(k). \label{eq:OFGF:opt_prob}
\end{equation}
Since the loss is assumed bounded below and proper, the infimum $\underline{\mathcal{L}}:=\inf_{k\in\mathcal{X}}\mathcal{L}(k)$ is attained; let $\mathcal{T}:=\{k\in\mathcal{X}~|~\mathcal{L}(k)=\underline{\mathcal{L}}\}$ be the \emph{target set}, or the set of points that solve \eqref{eq:OFGF:opt_prob}.
Define the optimality gap $\theta(k) \coloneqq \mathcal{L}(k)-\underline{\mathcal{L}} \ge 0$ and consider the gradient flow
\begin{equation}\label{eq:GF}
    \dot{k}(t)=-\nabla \mathcal{L}(k(t)),\qquad k(0)\in\mathcal{X},
\end{equation}
for finding a solution of \eqref{eq:OFGF:opt_prob}. One can verify that under a gradient flow it holds that $\dot \theta(t)=-\|\nabla \mathcal{L}(k(t))\|^2$, and since $\mathcal{L}$ is smooth on the open set $\mathcal{X}$, every global minimizer is stationary, hence
$\mathcal{T}\subseteq\mathcal{Z}:=\{k\in\mathcal{X}:\nabla\mathcal{L}(k)=0\}$.

Despite that, notice that in general stationarity is only necessary for optimality of a point, and convergence of a gradient method to global minima requires additional structure. Ideally one would characterize the following property for the gradient flow \eqref{eq:GF}:

\begin{definition}[global exponential cost stability (GECS)]\label{def:GES}
    Let $\phi:\mathbb{R}_+\times \mathcal{X}\to\mathcal{X}$ be the solution of \eqref{eq:GF} for any initial condition in $\mathcal{X}$. The  gradient flow \eqref{eq:GF} is said to be \emph{globally exponentially cost stable} (GECS) if $\exists\mu>0$ for which
    \begin{equation}\label{eq:GES}
        \theta\big(\phi(t,k_0)\big) \le \theta(k_0)\,e^{-\mu t}, \qquad \forall t\ge 0,\ \forall k_0\in\mathcal{X}.
    \end{equation}
\end{definition}

The following definition and lemma introduce necessary and sufficient conditions for this property to hold.

\begin{definition}[global Polyak--\L{}ojasiewicz (\gPLI)]\label{def:gPLI}
    A function $\mathcal{L}$ satisfies a \gPLI if there exists a $\mu > 0$ for which
    \begin{equation}\label{eq:gPLI}
        \|\nabla \mathcal{L}(k)\| \ge \sqrt{\mu (\mathcal{L}(k)-\underline{\mathcal{L}})}, \qquad \forall k\in\mathcal{X}.
    \end{equation}
\end{definition}

\begin{lemma}[Lemma 2 of \cite{ArthurCDC25}]\label{lem:gPLI-GES}
    The gradient flow \eqref{eq:GF} of a minimization problem \eqref{eq:OFGF:opt_prob} is GECS if and only if its loss function satisfies a \gPLI.
\end{lemma}

Traditionally, the \gPLI (or just P\L{}I as it is usually referred to in the literature \cite{polyak1963gradient,karimi2016linear}) is a powerful tool in nonconvex optimization, linking the gradient norm $\|\nabla \mathcal{L}(k)\|$ to the gap $\theta(k)$ and thus quantifying exponential decay rates for solutions of \eqref{eq:GF}. The \gPLI, however, can be impossible to characterize for some problems, and weaker versions of the inequality could be enough to provide weaker convergence guarantees. An example is the saturated P\L{}I condition, defined as follows:

\begin{definition}[Saturated P\L{}I (\satPLI)]\label{def:satPLI}
    A function $\mathcal{L}$ satisfies \emph{\satPLI} if there exist $a,b>0$ for which
    \begin{equation}\label{eq:satPLI}
        \|\nabla \mathcal{L}(k)\|\ \ge\ \sqrt{\frac{a(\mathcal{L}(k)-\underline{\mathcal{L}})}{b+(\mathcal{L}(k)-\underline{\mathcal{L}})}}, \qquad \forall k\in\mathcal{X}.
    \end{equation}
\end{definition}

In \cite{cui2024small} the authors showed that a policy-optimization strategy for the continuous-time linear quadratic regulator satisfies a \satPLI as defined. Furthermore, in \cite{ArthurCDC25}, we showed that the same problem will never satisfy a \gPLI.

A priori, if the cost satisfies a \satPLI, then only asymptotic convergence can be guaranteed. However, two asymptotically convergent solutions can differ greatly in qualitatively behavior, and a specific type of convergence behavior can be characterized for cost functions that satisfy a \satPLI.

\begin{definition}[global linear–exponential cost stable (GLECS)]\label{def:GLES}
    Let $\phi:\mathbb{R}_+\times \mathcal{X}\to\mathcal{X}$ be the solution of \eqref{eq:GF} for any initial condition in $\mathcal{X}$. The gradient flow \eqref{eq:GF} is \emph{GLECS} if $\exists\beta>0$ such that, for every $k_0\in\mathcal{X}$, and every $t^\ast\!\ge 0$ there exists $\mu_{k_0,t^\ast}\!>\!0$ for which
    \begin{equation}\label{eq:GLES}
      \theta(\phi(t,k_0)) \le
      \begin{cases}
        \theta(\phi(t^\ast,k_0))-\beta (t - t^\ast), & t\le t^\ast,\\[2pt]
        \theta(\phi(t^\ast,k_0))\,e^{-\mu_{k_0,t^\ast}(t-t^\ast)}, & t>t^\ast.
      \end{cases}
    \end{equation}
\end{definition}
\begin{figure}
    \centering
    \includegraphics[width=0.5\linewidth]{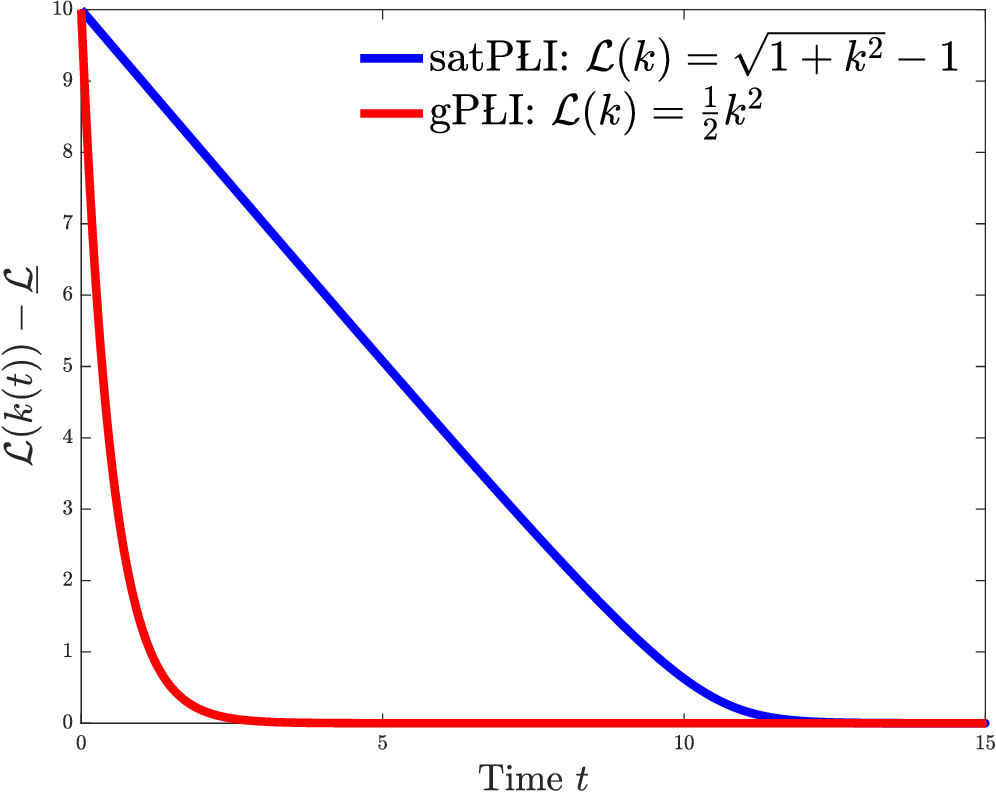}
    \caption{Depiction of two gradient flow solutions, one whose cost satisfies a \gPLI, and one whose cost satisfies only a \satPLI and with a globally bounded gradient.}
    \label{fig:GLECSvsGECS}
\end{figure}
\begin{lemma}[Lemma 3 of \cite{ArthurCDC25}]\label{lem:satPLI-GLES}
    If the gradient flow \eqref{eq:GF} of a minimization problem \eqref{eq:OFGF:opt_prob} satisfies a \satPLI \eqref{eq:satPLI} with some $a,b>0$ and there exists $L_\infty>0$ for which $\|\nabla \mathcal{L}(k)\|\le L_\infty$ holds for every $k\in\mathcal{X}$, then the flow is GLECS.
\end{lemma}



Intuitively, the solution of a GLECS gradient flow has linear convergence away from the optimal cost, and then does a ``soft-switch'' to exponential once it gets near the optimal value. Fig. \ref{fig:GLECSvsGECS} illustrates this convergence profile, and compares two different solutions, one whose cost satisfies a \gPLI and one that satisfies a \satPLI.

It is easy to verify, as stated in \cite{ArthurCDC25}, that a \gPLI implies a \satPLI but the converse is not true. In this paper, we bring attention to the fact that, through an appropriate reparameterization, a \satPLI problem can be reformulated into one satisfying a \gPLI. To illustrate this, we next present a simple case where a seemingly trivial reparameterization converts a problem whose gradient flow originally only satisfied a \satPLI into one that can be shown to satisfy a \gPLI.

\subsection{Motivating example}\label{subsec:Lift-sat-to-g}


Let $\mathcal{X}=(h,\infty)$, $h>0$. Assume $\mathcal{L}:\mathcal{X}\to\mathbb{R}$ is real-analytic, bounded below, proper, satisfies a \satPLI with constants $a,b>0$, and is such that $\limsup_{k\to\infty}\mathcal{L}'(k)$ is finite. For this function, let $k^\star$ be the largest element of $\mathcal{T}:=\{k\in(h,\infty)~|~\mathcal{L}(k)=\underline{\mathcal{L}}\}$, and define $L_\infty:=\max_{k\in[k^\star,\infty)}\mathcal{L}'(k)$. It is easy to see from the results presented previously that a solution of a gradient flow for this cost function initialized at points larger than $k^\star$ would present a linear-exponential profile similar to the one described in Definition \ref{def:GLES}. 

Now, consider the following reparameterization
\begin{equation*}
    f(k)\coloneqq \mathcal{L}(k^2),\qquad k\ge \sqrt{h},
\end{equation*}
which implies that $\underline f:=\inf_{k\ge \sqrt h} f(k)=\underline{\mathcal L}$ and $f'(k)=2k\,\mathcal L'(k^2)$. It is easy to see that both optimization problems are trivially equivalent, i.e., a point $k^*$ minimizes $f$ if and only if $(k^*)^2$ minimizes $\mathcal{L}$, with both cost functions attaining the same minimum value. 
Nonetheless, the gradient flow for $f$ admits the following result:



\begin{theorem}\label{thm:global-mu-eps}
    For any $\epsilon>\sqrt{h}$, there exists a $\mu_\epsilon>0$ such that for all $k\in[\epsilon,\infty)$, the reparameterized loss function $f$ satisfies the following \gPLI
    \begin{equation}\label{eq:gPLI-f}
        |f'(k)| \ge \sqrt{\mu_\epsilon(f(k)-\underline f)}\,.
    \end{equation}
\end{theorem}

Furthermore, the condition can be shown to hold in $(\sqrt{h},\infty)$ if further assumptions on $\mathcal{L}$ are made to guarantee the existence of the limit $\lim_{k\to\sqrt{h}}\mathcal{L}'(k)$. By Lemma~\ref{lem:gPLI-GES}, \eqref{eq:gPLI-f} yields a global exponential decay for the gradient flow of $f$, whereas the original $\mathcal{L}$ under \satPLI exhibits the linear–exponential behavior of Lemma~\ref{lem:satPLI-GLES} for solutions initialized at a value larger than $k^\star$. Thus, the reparameterization $k\mapsto k^2$ can qualitatively change the convergence behavior of solutions of a gradient flow. To further explore this phenomenon, we consider a simple and algebraically trackable case of the policy optimization LQR for continuous time system, whose standard formulation is known to satisfy only a \satPLI and never a \gPLI \cite{ArthurCDC25}. 

\section{The Overparameterized LQR Problem}\label{sec:LQRPolicyOptimization}

This section examines how overparameterization reshapes the convergence of gradient flow solutions. To better understand the behavior of policy optimization in this overparameterized LQR, we study a scalar system setup shown in Problem~\ref{gPLI:problem2} that allows for tractable analysis while preserving the essential dynamics of the general case. 

\begin{problem}\label{gPLI:problem2}
    Consider the scalar LTI system:
    \begin{equation}\label{eq:vec:dyn}
        \dot{x} = ax + u, \quad x(0) = 1,
    \end{equation}
    with arbitrary system parameter $a \in \mathbb{R}$. Consider an overparameterized state feedback, defined as $u = -\mathbf{k}x$ where $\mathbf{k} := k_2k_1\in\mathbb R$ and $k_1, k_2^\top \in \mathbb{R}^\kappa$, for some integer $\kappa>0$. The admissible set of control parameters $(k_1, k_2)$ is defined as $\mathcal{K} := \left\{ (k_1, k_2) \in \mathbb{R}^{\kappa \times 1} \times \mathbb{R}^{1 \times \kappa} \ \middle| \ a - \mathbf{k} < 0 \right\}$ to ensure closed-loop stability. We define the overparameterized LQR cost $\mathcal{L}:\mathcal{K}\to\mathbb{R}$ as $\mathcal{L}(k_1,k_2) \coloneqq J(\mathbf{k})$ with
    \begin{align}\label{eq:vec:cost-J}
        J(\mathbf{k}) &= \mathbb{E}_{x_0\sim\mathcal{X}} \left[ \int_0^\infty \left( x^2 q + u^2 r \right) dt \right] \nonumber\\
        &= \mathbb{E}[x_0^2] \cdot \left[q + \mathbf{k}^2 r \right] \int_0^\infty e^{2(a - \mathbf{k})t} dt \nonumber\\
        &= -\frac{q + r\mathbf{k}^2}{2(a - \mathbf{k})} > 0,
    \end{align}
    where $\mathcal{X}\sim\mathcal{N}(0,1)$. Notice that convergence requires \( \mathbf{k} > a \), and the cost weights to satisfy $q,r > 0$.
\end{problem}
For this simple version of the problem, we can compute the gradient of \( J(\mathbf{k}) \) with respect to $\mathbf{k}$ explicitly as
\begin{equation}\label{eq:vec:grad-J}
    \nabla J(\mathbf{k}) := \frac{r\mathbf{k}^2 - 2ar\mathbf{k} - q}{2(a - \mathbf{k})^2},
\end{equation}
and using the chain rule we obtain
\begin{equation*}
    \nabla_{k_1} \mathcal{L}(k_1,k_2) = \nabla J(\mathbf{k}) k_2^\top, \quad \nabla_{k_2} \mathcal{L}(k_1,k_2) = \nabla J(\mathbf{k}) k_1^\top.
\end{equation*}
Notice that the evolution of the cost along trajectories of \( (k_1, k_2) \) will satisfy
\begin{equation}\label{eq:vec:grad-L}
    \dot{\mathcal{L}} \coloneqq - \lVert \nabla \mathcal{L}(k_1,k_2)\rVert^2 = -\nabla J(\mathbf{k})^2(\|k_1\|^2 + \|k_2\|^2),
\end{equation}
which is not enough to guarantee global convergence to the optimal solution. To guarantee optimality of the solution, we aim to find the largest subset of $\mathcal{K}$ in which we can establish a global Polyak–\L{}ojasiewicz inequality (\gPLI), i.e.
\begin{equation}\label{eq:vec:gPLI}
    \|\nabla \mathcal{L}(k_1,k_2)\| \geq \sqrt{\mu(\mathcal{L}(k_1,k_2) - \underline{\mathcal{L}})},
\end{equation}
for all \( (k_1,k_2) \in \overline{\mathcal{K}} \), where \( \underline{\mathcal{L}} := \inf_{(k_1,k_2) \in \mathcal{K}} \mathcal{L}(k_1,k_2) \) and $\overline{\mathcal{K}}\subseteq\mathcal{K}$. 


To achieve this goal, we first define a known property of overparameterized formulations \cite{ArthurAutomatica25,min2023convergence,menon2024geometry},

\begin{definition}[Imbalance Measure]\label{def:imbalance}
    For the overparameterized LQR in Problem~\ref{gPLI:problem2}, the invariant matrix $\mathcal{C}$ is defined as $\mathcal{C} := k_1 k_1^\top - k_2^\top k_2 \in \mathbb{R}^{\kappa \times \kappa}$, which remains constant along any solution of the gradient flow. Furthermore, the \emph{imbalance measure} is defined as $c := 2\,\tr[\mathcal{C}^2] - \left(\tr[\mathcal{C}]\right)^2$.
\end{definition}

The invariant $\mathcal{C}$ is called this because it is invariant along any solution of the overparameterized gradient flow. This is a well known property of overparameterization, and is proven for the LQR in \cite{ArthurAutomatica25}. The imbalance measure $c$ quantifies the asymmetry: it vanishes when the two factors $(k_1,k_2)$ are perfectly aligned (up to sign) and increases as asymmetry grows. Notice that for the scalar case, the invariant $\mathcal{C}$ is the difference of two rank-one, positive semi-definite matrices, for which the following result can be stated.

\begin{proposition}\label{prop:imbalance}
    Let $A,B\in\mathbb{R}^{n\times n}$ be any two rank-one, positive semi-definite matrices, and let $C=A-B$. Then $c:=2\tr[C^2]-(\tr[C])^2\ge0$, with $c=0$ if and only if $A=B$.
\end{proposition}

Proposition~\ref{prop:imbalance} guarantees that the imbalance measure \( c \) is always nonnegative. Furthermore, $c=0$ if and only if $k_1=\pm k_2$, which means that the center-stable manifold of the saddle point lies entirely within the set $\{(k_1,k_2)\in\mathcal{K}~|~c(k_1,k_2) = 0\}$.
With this in mind, and with the the constraint set \( \mathcal{K}_\gamma \) we can state the following result

\begin{theorem}\label{thm:gPLI}
    Consider the overparameterized model-free LQR in Problem~\ref{gPLI:problem2} and for any $\gamma>0$, define the set $\mathcal{K}_\gamma \coloneqq \left\{ (k_1, k_2)\in\mathcal{K} \,\middle|\, d(k_1,k_2)\coloneqq\|k_1 + k_2^\top\|^2 \geq \gamma \right\}$. Then, for \( \gamma > \max(0,4a)\):
    \begin{enumerate}
        \item any solution initialized in $\mathcal{K}_\gamma$ remains in $\mathcal{K}_\gamma$; and
        \item the \emph{\gPLI} $\|\nabla \mathcal{L}(k_1,k_2)\| \geq \sqrt{\mu_\gamma(\mathcal{L}(k_1,k_2) - \underline{\mathcal{L}})}$ holds for all $(k_1,k_2)\in\mathcal{K}_\gamma$ with
        \begin{equation*}
            \mu_\gamma(c, \gamma) = \begin{cases}
                \frac{r}{4} \sqrt{ \frac{4c + \gamma^2}{(\gamma - 4a)^2} } > 0 & \textrm{ if } a \ge 0, \\[.5em]
                \frac{r}{4} \sqrt{ \frac{(\gamma^2 + c)^2}{(\gamma^2 - c - 4a\gamma)^2} } > 0 & \textrm{ if } a < 0, c < \Tilde{c} \\[.5em]
                \frac{r}{4} \sqrt{ \frac{c}{4a^2 + c} } > 0 & \textrm{ if } a < 0, c \ge \Tilde{c}
            \end{cases}
        \end{equation*}
        where \( \Tilde{c} \coloneqq \frac{a\gamma^2}{a - \gamma} > 0 \).
    \end{enumerate}
\end{theorem}


From Theorem~\ref{thm:gPLI}, the function $\mu_\gamma(c,\gamma)$ depends jointly on the imbalance measure $c$ and the set $\mathcal{K}_\gamma$ with parameter $\gamma$, for arbitrary $a$. Thus, when we seek a \emph{uniform guarantee} that the \gPLI holds for all admissible $(k_1,k_2) \in \mathcal{K}_\gamma$, we take the smallest possible value of $\mu_\gamma(c,\gamma)$ across all $c \ge 0$ for a fixed $\gamma$, yielding the expression in Corollary~\ref{cor:mu-gPLI}.

\begin{corollary}\label{cor:mu-gPLI}
    Consider the overparameterized model-free LQR in Problem~\ref{gPLI:problem2} and the set $\mathcal{K}_\gamma \coloneqq \left\{ (k_1, k_2)\in\mathcal{K} \,\middle|\, d\geq \gamma \right\}$. Then, for all \( \gamma > \max(0,4a) \), there exists a smallest constant 
    \begin{equation}\label{eq:vec:lowerbound}
        \underline{\mu} \coloneqq \frac{r}{4} \min\left\{1, \sqrt{\frac{\gamma^2}{(\gamma - 4a)^2}}\right\} > 0,
    \end{equation}
    such that the \emph{\gPLI} $\|\nabla \mathcal{L}(k_1,k_2)\| \geq \sqrt{\underline{\mu}(\mathcal{L}(k_1,k_2) - \underline{\mathcal{L}})}$ holds for all $(k_1,k_2)\in\mathcal{K}_\gamma$.
\end{corollary}


Furthermore, the function $\mu_\gamma(c,\gamma)$ in Theorem~\ref{thm:gPLI} prompts the question: \emph{how does the imbalance of the factorization $(k_1,k_2)$ affect the convergence rate within the admissible set $\mathcal{K}_\gamma$?} To answer this question, fix $\gamma$ and compare two trajectories with equal initial cost but different factorizations (hence different $c$). 
Intuitively, if $\mu_\gamma(c,\gamma)$ is a monotonic function in $c$ (for a fix $\gamma$), then a trajectory with greater imbalance $c$ will increases or decrease $\mu_\gamma(c,\gamma)$, and hence affecting the decay rate. Thus, the next corollary formalizes this intuition by analyzing the function $\mu_\gamma(c,\gamma)$.

\begin{corollary}\label{cor:convergence}
    Consider the overparameterized model-free LQR in Problem~\ref{gPLI:problem2}. Let \( (\phi_{k_1}(t), \phi_{k_2}(t)) \) be the solution of the gradient flow \eqref{eq:vec:grad-L} initialized at \( (k_1, k_2) \), and let \( \mathcal{L}(\phi_{k_1}(t), \phi_{k_2}(t)) \) be the corresponding cost trajectory. Suppose two initializations \( (\bar{k}_1, \bar{k}_2), (\hat{k}_1, \hat{k}_2) \in \mathcal{K}_\gamma \) have imbalance measures \( \bar{c} \) and \( \hat{c} \), respectively, satisfying $\mathcal{L}(\bar{k}_1, \bar{k}_2) = \mathcal{L}(\hat{k}_1, \hat{k}_2)$ and $ \bar{c} > \hat{c}$. Then the associated functions satisfy \( \mu_{\gamma}(\bar{c},\gamma) > \mu_{\gamma}(\hat{c},\gamma) \), and the trajectories obey
    \begin{equation*}
        e^{-\mu_{\gamma}(\bar{c},\gamma) t} \big(\mathcal{L}(\bar{k}_1, \bar{k}_2) - \underline{\mathcal{L}}\big) < e^{-\mu_{\gamma}(\hat{c},\gamma) t} \big(\mathcal{L}(\hat{k}_1, \hat{k}_2) - \underline{\mathcal{L}}\big),
    \end{equation*}
    showing that greater imbalance leads to faster convergence bounds.
\end{corollary}


The next section presents numerical experiments indicating that the formulation-dependent qualitative change from GLECS to GECS observed in the scalar model also persists for general LQR instances. 

\section{Numerical Results}\label{sec:NumericalResults}

In this section we present numerical simulations to demonstrate the theoretical results in a more general setting. Our aims are twofold:  
(i) to demonstrate the \emph{qualitative} change proven for the scalar case in Corollary \ref{cor:mu-gPLI}; and  
(ii) illustrate the trade-off between formulations generated by the interplay of the previous result and the introduction of saddle-points by the overparameterized formulation.

We consider two continuous-time LTI systems, one with stable system $\mathcal{G}_1$ and one with unstable system $\mathcal{G}_2$, given by 
\begin{equation*}
    \mathcal{G}_1 \coloneqq \{\dot{x}_1 = A^- x_1 + B_1 u_1\}, \quad \mathcal{G}_2 \coloneqq \{\dot{x}_2 = A^+ x_2 + B_2 u_2\},
\end{equation*}
where $x_1,x_2\in\mathbb{R}^{n}$ is the state with the initial conditions $x_{1,2}(0)\sim\mathcal{N}(0,I_n)$ and $I_n$ is the identity matrix of size $n$. Both $u_1\in\mathbb{R}^{m_1}$ and $u_2\in\mathbb{R}^{m_2}$ are the control inputs given by the overparameterized neural network. Both inputs are denoted by $u = \mathbf{K}x$ (LFFNN), with $\mathbf{K} = K_2K_1$, $K_1\in\mathbb{R}^{\kappa\times n}$, $K_2\in\mathbb{R}^{m\times \kappa}$ where $m$ matches the number of inputs, either $m_1$ or $m_2$, and $\kappa = 10$ defines the number of neurons. The matrices are defined as $A^- \coloneqq -(5I_n + A)\in\mathbb{R}^{n\times n}$ and  $A^+ \coloneqq A\in\mathbb{R}^{n\times n}$, with
\begin{equation*}
    A = 
    \begin{bmatrix}
    0.2373 & 0.3452 & 0.6653 & 0.6715 & 0.3288 \\
    0.3452 & 0.4889 & 0.8060 & 0.3889 & 0.5584 \\
    0.6653 & 0.8060 & 0.0377 & 0.5735 & 0.5100 \\
    0.6715 & 0.3889 & 0.5735 & 0.3354 & 0.6667 \\
    0.3288 & 0.5584 & 0.5100 & 0.6667 & 0.4942
    \end{bmatrix},
\end{equation*}
while $B_1\in\mathbb{R}^{n\times m_1}$ and $B_2\in\mathbb{R}^{n\times m_2}$, denoted as
\begin{equation*}
    B_1^\top =
    \begin{bmatrix}
    0 & 0 & 1 & 0 & 0 \\
    0 & 0 & 0 & 1 & 0 \\
    0 & 0 & 0 & 0 & 1
    \end{bmatrix},
    \quad
    B_2 = I_n,
\end{equation*}
and $m_2 = n$. 

We consider the overparameterized LQR problem as defined in \cite{ArthurAutomatica25}, with $J(\mathbf{K})$ being the LQR cost and $\mathcal{L}(K_1,K_2):=J(\mathbf{K})$ being its overparameterized formulation.
%
The evolution of $\mathcal{L}$ along $(K_1,K_2)$ is defined by the following dynamics
\begin{equation*}
    \dot{K}_1 = -K_2^\top \nabla J(\mathbf{K}), \quad \dot{K}_2 = -\nabla J(\mathbf{K})K_1^\top.
\end{equation*}

For both systems, $\mathcal{G}_1$ and $\mathcal{G}_2$, we take an initial feedback matrix $\mathbf{K}^-(0) \in \mathbb{R}^{m_1\times n}$ and $\mathbf{K}^+(0) \in \mathbb{R}^{m_2\times n}$, each guaranteeing that their respective close-loop system are stable. 
In the Ovp. LQR setting, these initial gains are factorized as $\mathbf{K}(0) = K_2(0)K_1(0)$, with both $K_1(0)$ and $K_2(0)$ chosen to match the initialization of each scenario, $\mathcal{G}_1$ and $\mathcal{G}_2$ (see \cite{ArthurAutomatica25} or  Remark~\ref{rem:factorization} as an example on how to factorize $\mathbf{K}(0)$). The two initializations are
\begin{equation*}
    \mathbf{K}^{-}(0) = \begin{bmatrix}
        -2.14  &  -2.62  &  20.48  &  -1.55  &  -1.30 \\
        -2.07  &  -0.80  &  -1.55  &  19.14  &  -1.94 \\
        -0.64  &  -1.50  &  -1.30  &  -1.94  &  18.44
    \end{bmatrix}
\end{equation*}
\begin{equation*}
    \mathbf{K}^{+}(0) = \begin{bmatrix}
    12.21  &  1.12  &  2.16  &  2.18  &  1.07 \\
    1.12  &  13.03  &  2.61  &  1.26  &  1.81 \\
    2.16  &  2.61  &  11.56  &  1.86  &  1.65 \\
    2.18  &  1.26  &  1.86  &  12.53  &  2.16 \\
    1.07  &  1.81  &  1.65  &  2.16  &  13.02
    \end{bmatrix}.
\end{equation*}
The optimal value of the cost is denoted by $\underline{\mathcal{L}}$ while the optimality gap is defined by $\mathcal{L}(\mathbf{K}(t)) - \underline{\mathcal{L}}$ for the Ovp. LQR and $\mathcal{L}(K(t)) - \underline{\mathcal{L}}$ for standard LQR.

\begin{figure}[t!]
    \centering
    \subfloat[\label{Fig:NumRes1b} Stable case $\mathcal{G}_1$]{\includegraphics[width=.475\linewidth]{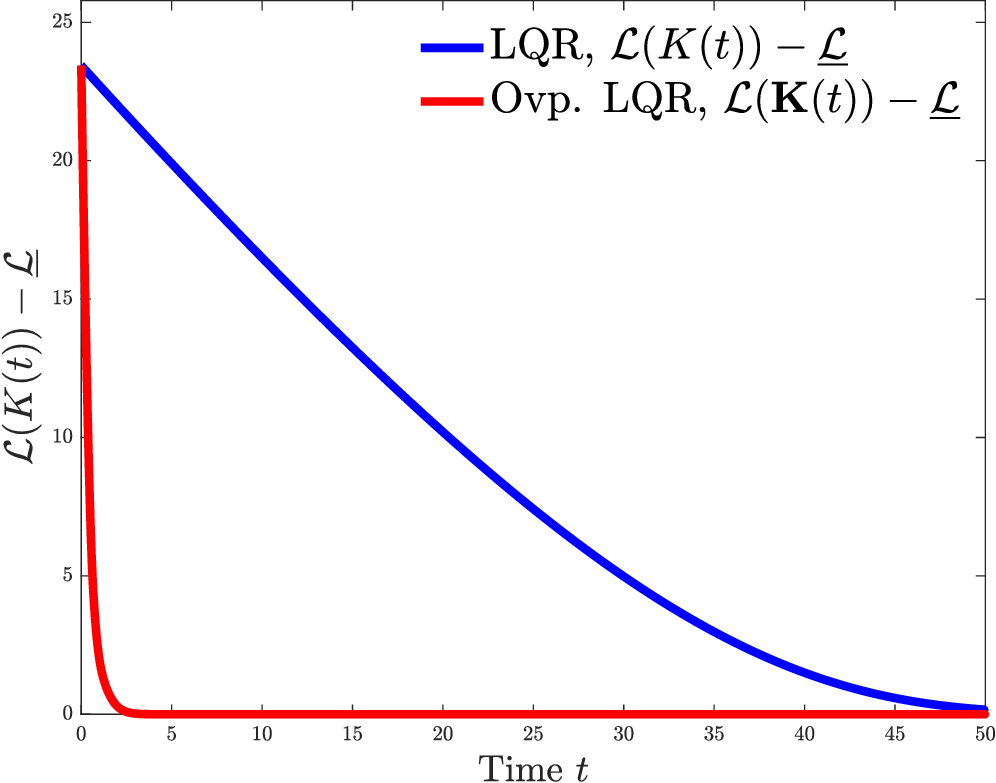}}\,
    \subfloat[\label{Fig:NumRes1d} Unstable case $\mathcal{G}_2$]{\includegraphics[width=.475\linewidth]{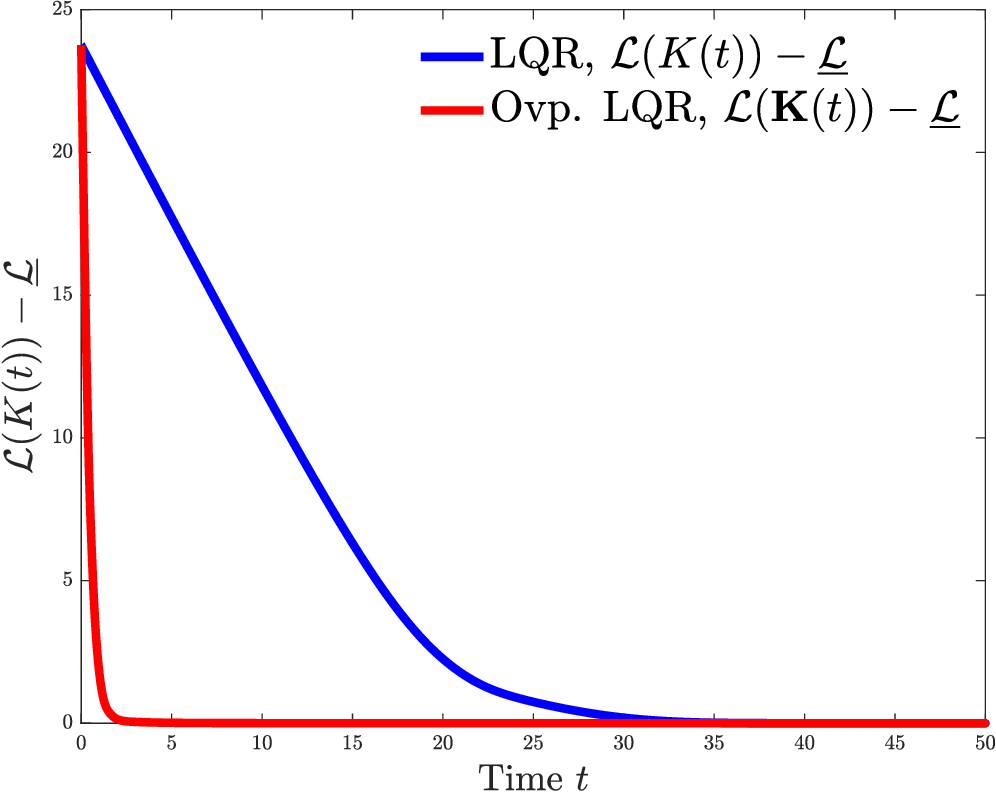}}
    \caption{Optimality–gap comparison for standard LQR vs. Ovp. LQR on two systems, $\mathcal{G}_1$ and $\mathcal{G}_2$. Both methods start with the same initial cost $\mathbf{K}^-(0) \!=\! K^-(0)$ and $\mathbf{K}^+(0) \!=\! K^+(0)$ with initial states $x_1(0)$ and $x_2(0)$ in turn.}
    \label{Fig:NumRes1}
\end{figure}

We discuss the first point: the Ovp.\ LQR exhibits GECS, whereas the standard LQR shows a GLECS. This is portrayed in Figs.~\ref{Fig:NumRes1b}–\ref{Fig:NumRes1d}, comparing the optimality gap of the standard LQR policy optimization (blue) against the Ovp. LQR counterpart (red) for the stable $\mathcal{G}_1$ and the unstable system $\mathcal{G}_2$. Consistent with Theorem~\ref{thm:gPLI} and Corollary~\ref{cor:mu-gPLI}, the Ovp.\ LQR exhibits an exponential decay of the cost gap (GECS), while the standard LQR shows the characteristic linear--exponential profile (GLECS). The effect is pronounced in both the stable and unstable regimes, underscoring that the qualitative improvement is due to the reparameterization rather than system stability. Empirically, the factorization $(K_1,K_2)$ reshapes the geometry so that the policy gradient through the chain rule can be larger in descent directions of $K$, yielding a faster convergence decay.

\begin{figure}[t!]
    \centering
    \includegraphics[width=0.475\linewidth]{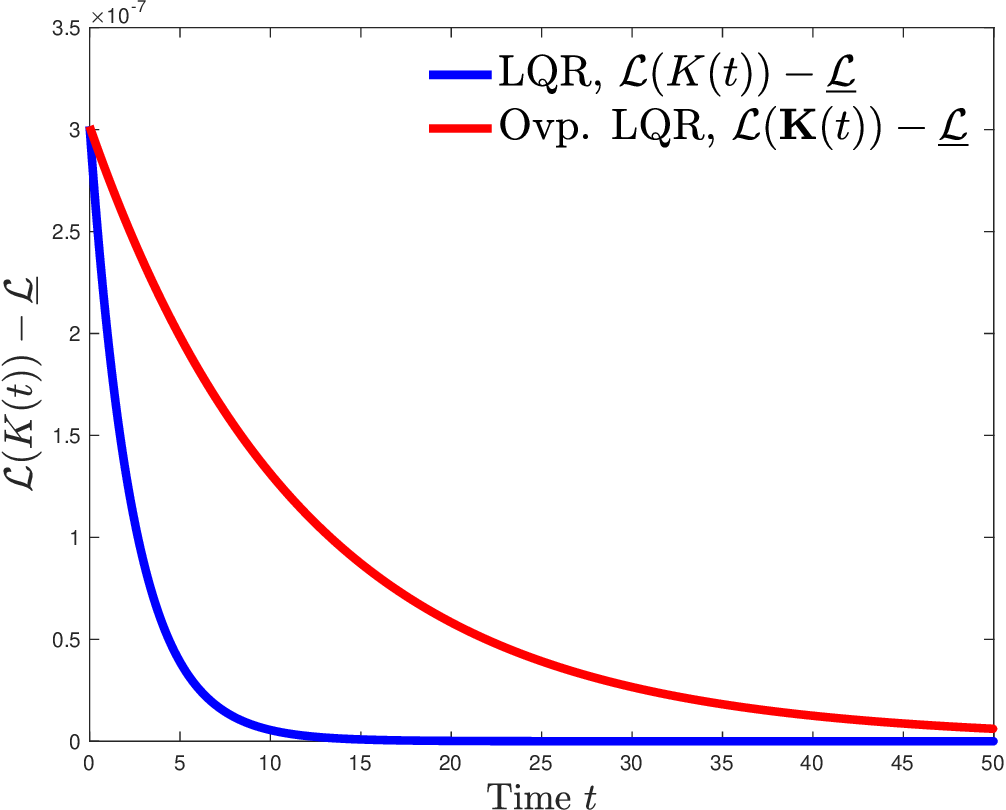}
    \caption{The optimality–gap in $\mathcal{G}_1$ using $\mathbf{K}^-(0) = K^-(0) \approx 0_{m,n}$}
    \label{fig:saddle}
\end{figure}

Figure~\ref{fig:saddle} shows the second point using $\mathcal{G}_1$: what happens when the Ovp.\ initialization is close to the center–stable manifold of the saddle $(\gamma\approx 0$, hence $c \approx 0)$, in which we pick $\mathbf{K}(0)\approx 0$. According to Theorem~\ref{thm:gPLI} and Corollary~\ref{cor:mu-gPLI}, this results in a small value of the function $\mu_\gamma(c,\gamma)$, yielding slow convergence (often slower than the standard LQR) for the gradient flow but consistent with exponential decay. Practically, this highlights a trade–off: overparameterization yields GECS, but care is needed to avoid the saddle (e.g., by enforcing $d\ge\gamma$ at initialization) to prevent slow convergence.

\begin{remark}[Initialization of $K_1(0)$ and $K_2(0)$]\label{rem:factorization}
    Define the initial cost gap $\eta>0$. Let \(K^\ast\) be the LQR gain for \((A,B,Q,R)\). Choose a scale \(s>0\) and set \(\mathbf{K}^\ast:=sK^\ast\). Accept \(\mathbf{K}^\ast\) only if \(A_{\mathrm{cl}}:=A-B\mathbf{K}^\ast\) is Hurwitz, in which case the cost $\mathcal{L}(\mathbf{K}^\ast)=\operatorname{trace}(P_0)$ where $P_0$ solves
    \begin{equation*}
        A_{\mathrm{cl}}^{\top} P_0 + P_0 A_{\mathrm{cl}}+\bigl(Q +\mathbf{K}^{\ast\top}R\mathbf{K}^\ast\bigr) = 0,
    \end{equation*}
    yields the gap \(\theta_0 \coloneqq \mathcal{L}(\mathbf{K}^\ast)-\underline{\mathcal{L}}\). Increase \(s\) geometrically until the gap exceeds the target $\theta_0 \ge \eta$, then set \(\mathbf{K}(0):=\mathbf{K}^\ast\).
    
    Draw any a full–column–rank matrix \(K_1(0)\in\mathbb{R}^{\kappa\times n}\) with \(\kappa\ge n\) and define $K_2(0):=\mathbf{K}(0)(K_1(0)^{\!\top}K_1(0))^{-1}K_1(0)^{\!\top}$. Hence \(K_2(0)K_1(0)=\mathbf{K}(0)\). This initialization (i) ensures a controlled initial gap relative to the optimum and (ii) matches the desired overparameterized product exactly.
\end{remark}

\section{Conclusion and Future Work}\label{sec:Conclusion-FutureWork}
We proved that problem formulation shapes the convergence behavior of gradient flows in nonconvex optimization. A simple reparameterization $k\mapsto k^2$ converts a landscape that only satisfies only a \satPLI into one that satisfies a \gPLI \emph{almost} everywhere, showing a \emph{qualitative} change in convergence for the solution; from global linear–exponential cost stability (GLECS) to global exponential cost stability (GECS). Motivated by this observation, we consider the policy-optimization LQR as a testbed, and showed that an overparameterized (LFFNN) policy achieves the same \emph{qualitative} upgrade. For the scalar testbed, we characterized a \gPLI with an \emph{explicit} function rate that depends on an ``imbalance'' measure $c\ge 0$ and an invariant region $\mathcal{K}_\gamma$ in which solutions must be initialized. Numerical experiments in a general setting on stable and unstable systems confirmed the theory: the overparameterized formulation exhibits GECS while the standard formulation shows GLECS; and highlighted the expected slow convergence only when initialized near the saddle manifold set $c \approx 0$. Future works will extend these guarantees to stochastic gradients and noisy dynamics; generalize beyond LQR to broader control/learning problems; and analyze deeper/alternative factorizations to sharpen \gPLI constants and robustness.

\bibliographystyle{ieeetr}  
\bibliography{Ref} 

\appendix

\section{Proofs}
    \subsection*{Proof of Theorem \ref{thm:global-mu-eps}}
    This proof will be derived for $k$ in the domain of $f$, not $\mathcal{L}$, as such, all optimal solutions and boundaries are such for $f$ and their square are the equivalent for $\mathcal{L}$. Furthermore, we assume for simplicity that $\underline{f}=\underline{\mathcal{L}}=0$, but the results here can be adapted for arbitrary minimum values. 
    
    We divide the proof into two parts, one for $k\in(\epsilon,k_*)$ and one for $k\in(k^*,\infty)$ where $k_*$ is the smallest element of $\mathcal{T}:=\{k\in(\epsilon,\infty)~|~f(k)=0\}$, and $k^*$ is the largest. We show that for each case, a $\mu>0$ exists, and so a $\mu$ that holds everywhere can be found by taking the smallest of the two.

    \textbf{First case $k\in(\epsilon,k_*)$}: We want to show that there exists a $\rho_*>0$ such that
    \begin{equation*}
        \psi(k):=\frac{(f'(k))^2}{f(k)}\geq\rho_*
    \end{equation*}
    for all $k\in(\epsilon,k_*)$. We argue that fact through properness of $f$. Because $f$ is a proper function on the interval $(\sqrt{h},\infty)$, then it is strictly (because of the \satPLI) decreasing on the interval $(\sqrt{h},k_*)$ from $\infty$ to $0$. Notice that the function $\psi(k)$ can never be zero for any point on the interior of $(\sqrt{h},k_*)$, otherwise strict monotonicity of $f$ would not be true.
    
    Next, notice that because $\mathcal{L}(\cdot)$ satisfies a \satPLI and $f(k)=\mathcal{L}(k^2)$, we can write
    \begin{align*}
        \mu(k) &:=\frac{(f'(k))^2}{f(k)}=4k^2\frac{(\mathcal{L}'(k^2))^2}{\mathcal{L}(k^2)}\geq\frac{4ak^2}{b+f(k)}.
    \end{align*}
    which implies that implies that, $\lim_{k\to k_*}r(k)\geq{4ak_*^2}/{b}$. This proves that $\rho_*=\inf_{k\in[\epsilon,k_*)}\psi(k)>0$ always exists for all $\epsilon>\sqrt{h}$. 

    \textit{Remark:} We point here, without proof, that the result above can be generalized to the entire interval $(\sqrt{h},k_*)$, if the limit $\lim_{k\to\sqrt{h}}\psi(k)$ is assumed to exist.

    \textbf{Second case $k\in(k^*,\infty)$}: Since $h>0$, then $k^*>0$. From this and the asymmetrical gradient bound for $\mathcal{L}$, we can conclude that it holds that $L_\infty k^2\geq f(k)=\mathcal{L}(k^2)$ for all $k\in(k^*,\infty)$. Then notice that:
    \begin{align*}
        \frac{(f'(k))^2}{f(k)}&=4k^2\frac{(\mathcal{L}'(k^2))^2}{\mathcal{L}(k^2)}\geq 4k^2\frac{a}{b+f(k)}
    \end{align*}
    From here, we again divide the analysis in two. First, if $f(k)<b$ then it follows that
    \begin{align*}
        \frac{(f'(k))^2}{f(k)}&\geq 4k^2\frac{a}{b+f(k)}\geq\frac{2a(k^*)^2}{b}=:\rho_1
    \end{align*}
    and if $f(k)>b$, then
    \begin{align*}
        \frac{(f'(k))^2}{f(k)}&\geq 4k^2\frac{a}{b+f(k)}
        \geq\frac{4}{L_\infty}\frac{af(k)}{2f(k)}=\frac{2a}{L_\infty}=:\rho_2.
    \end{align*}
    Finally, notice that the global Polyak-\L{}ojasiewicz inequality holds for $k\in(k^*,\infty)$ for $\rho^*:=\min(\rho_1,\rho_2)>0$, and globally for $k\in(\sqrt{h},\infty)$ for $\rho=\min(\rho_*,\rho^*)>0$. \qed

    \subsection*{Proof of Proposition~\ref{prop:imbalance}}
    \newcommand{\aaa}{\mathbf{a}}
    \newcommand{\bb}{\mathbf{b}}
    Since $A,B$ are rank-one PSD, write $A=aa^\top$ and $B=bb^\top$ for some $a,b\in\mathbb{R}^n$. Let $\aaa:=a^\top a$ and $\bb:=b^\top b$ Then notice, after some algebraic manipulation, that
    \begin{align*}
        c \coloneqq 2\tr[C^2]-(\tr[C])^2 &= (\aaa+\bb)^2-4(a^\top b)^2 \\ &\geq(\aaa+\bb)^2-4\aaa\bb\\&=(\aaa-\bb)^2\geq 0,
    \end{align*}
    concluding the proof. \qed

\subsection*{Proof of Theorem~\ref{thm:gPLI}} 
    Before proceeding, to simplify the proofs, we introduce a scalar perturbation around the optimal feedback gain. Notice that the cost in \eqref{eq:vec:cost-J} and the gradient in \eqref{eq:vec:grad-J} can be rewritten as \( \Tilde{J}(\mathbf{k}) = p_\mathbf{k} \) and \( \nabla \Tilde{J}(\mathbf{k}) = -2(p_{\mathbf{k}} - r\mathbf{k})\ell \), where \( p_\mathbf{k} \) solves \( 2p_\mathbf{k}(a - \mathbf{k}) + r\mathbf{k}^2 + q = 0 \), and \( \ell \) satisfies \( 2\ell(a - \mathbf{k}) + 1 = 0 \). Let \( \mathbf{k}^\ast \) denote the optimal gain minimizing \( \Tilde{J}(\mathbf{k}) \) such that $\underline{\Tilde{J}} = \Tilde{J}(\mathbf{k}^\ast)$, and introduce a scalar perturbation \( \varepsilon \in \mathbb{R} \), so that \( \varepsilon := \mathbf{k} - \mathbf{k}^\ast\). Then, one can verify that:
    \begin{align}\label{eq:vec:varepsilon}\begin{aligned}
        \Tilde{J}(\mathbf{k}^\ast + \varepsilon) - \underline{\Tilde{J}} &= r \ell \varepsilon^2 \\
        \nabla \Tilde{J}(\mathbf{k}^\ast + \varepsilon) &= \ell(-r \ell \varepsilon^2 + r \varepsilon) \\
        f(\mathbf{k}^\ast + \varepsilon) &= r \ell (\ell^2 \varepsilon^2 - 2 \ell \varepsilon + 1)
        \end{aligned}
    \end{align}
    where $ f(\mathbf{k}^\ast + \varepsilon) \coloneqq \lVert \nabla \Tilde{J}(\mathbf{k}^\ast + \varepsilon) \rVert^2/[\Tilde{J}(\mathbf{k}^\ast + \varepsilon) - \underline{\Tilde{J}}]$ and $ \ell := -1/[2(a - \mathbf{k}^\ast - \varepsilon)] = -1/[2(a - \mathbf{k})]$. With these estabilished, we move on with the proofs. 
    
The first claim of the theorem follows from the geometry of the gradient flow and the role of $\mathcal{K}_\gamma$, and is proven in details in \cite{ArthurAutomatica25}. To prove the second claim, we aim to lower bound the constant $\mu_\gamma$ away from zero. To do so, notice that $\nabla_{\mathbf{k}} J(\mathbf{k}) = 0$ implies that \( \mathbf{k}^\ast = a + \sqrt{a^2 + \frac{q}{r}}\in\mathcal{K} \). Then, after some algebraic manipulation, we can write $\|k_1\|^2 + \|k_2\|^2 = 2\|k_2\|^2 + \tr[\mathcal{C}] = \sqrt{c + 4\mathbf{k}^2}$. Using that and \eqref{eq:vec:varepsilon}, we obtain:
\begin{equation}\label{eq:proof:thm:gPLI:function}
    \begin{aligned}
    \|\nabla \mathcal{L}(k_1,k_2)\|^2 
    &= \underbrace{\vartheta_1(\varepsilon)\vartheta_2(\mathbf{k})}_{\geq \mu_\gamma}\;(\mathcal{L}(k_1,k_2) - \underline{\mathcal{L}}) \end{aligned}
\end{equation}
where $\vartheta_1(\varepsilon)\coloneqq r(\ell^2\varepsilon^2 - 2\ell\varepsilon + 1)$ and $\vartheta_2(\mathbf{k}) := \ell \sqrt{c + 4\mathbf{k}^2}$. To prove the result, we then need to bound both $\vartheta_1$ and $\vartheta_2$ away from zero. We proceed with each function separately.

\textbf{Lower bound of \( \vartheta_1(\varepsilon) \) term}: Define \( \delta \coloneqq a - \mathbf{k}^\ast < 0 \), since \( \mathbf{k}^\ast > a \). From \eqref{eq:vec:varepsilon}, we also know that \( \ell = 1 / [2(\varepsilon - \delta)] > 0 \) since $\varepsilon = \mathbf{k} - \mathbf{k}^\ast$ and $\mathbf{k} > a$, yielding $\varepsilon > \delta$. Therefore, the product \( \bar{x} \coloneqq \ell \varepsilon = \varepsilon / [2(\varepsilon - \delta)] \) is well-defined for \( \varepsilon \in (\delta, \infty) \). As \( \varepsilon \to \delta^+ \), \( \bar{x} \to -\infty \), and as \( \varepsilon \to \infty \), \( \bar{x} \to \frac{1}{2} \). Furthermore, the derivative $\bar{x}'(\varepsilon) = \frac{-\delta}{2(\varepsilon - \delta)^2}$
is strictly positive for all \( \varepsilon > \delta \), since \( \delta < 0 \). Therefore, \( \bar{x} \) is strictly increasing over \( \varepsilon\in(\delta, \infty) \), and its image lies in \( (-\infty, \frac{1}{2}) \). Now observe that under coordinate change, $\vartheta_1(\bar{x}) = r(\bar{x} - 1)^2$. Since \( \bar{x} < \frac{1}{2} < 1 \), we know \( (\bar{x} - 1)^2 \) is strictly decreasing in \( \bar{x} \), and thus attains its minimum at \( \bar{x} = \frac{1}{2} \) (noting it is a parabola with minimum at \( \bar{x} = 1 \)). Therefore, this gives
\begin{equation}\label{eq:proof:thm:gPLI:v1}
    \vartheta_1(\varepsilon) \geq r\left( \frac{1}{2} - 1 \right)^2 = \frac{r}{4} > 0.    
\end{equation}

\textbf{Lower bound of the term} $\vartheta_2(\mathbf{k})$: To find the minimum of $\vartheta_2$, it suffices to minimize \( \vartheta_3(\mathbf{k}) \coloneqq (c + 4\mathbf{k}^2)\ell^2 \), since the square root is strictly increasing on $\mathbb{R}_+$. For $a\neq 0$, taking the derivative of \( \vartheta_3 \) with respect to \( \mathbf{k} \), we compute:
\begin{align*}
    \vartheta_3'(\mathbf{k})
    &= \frac{2\mathbf{k}}{(-a + \mathbf{k})^2} - \frac{2\left(\frac{c}{4} + \mathbf{k}^2\right)}{(-a + \mathbf{k})^3} = -\frac{4a\mathbf{k} + c}{2(-a + \mathbf{k})^3}.
\end{align*}
Setting $\vartheta_3'(\mathbf{k}_{\inf}) = 0$ gives $\mathbf{k}_{\inf}=-c/(4a)$, the unique critical point of \( \vartheta_3(\mathbf{k}) \). Whether this corresponds to a minimum of \( \vartheta_2 \) depends on the value of the parameters. 

For $a = 0$, $\vartheta_3(\mathbf{k})$ reduces to $\vartheta_3(\mathbf{k}) = \frac{c}{4\mathbf{k}^2} + 1$ meaning that \(\vartheta_3'(\mathbf{k}) = -\frac{c}{2\mathbf{k}^3}\leq0\) for all $\mathbf{k}>0$, since $c>0$, implying that the function has no critical points. Its values decrease from \(+\infty\) as \( \mathbf{k} \to 0^+ \) to \(1\) as \( \mathbf{k} \to \infty\). 


Next, we distinguish three cases based on the system parameter \( a \), as this determines the structure of the admissible domain \( \mathbf{k} \in (a, \infty) \) and whether the critical point \( \mathbf{k}_{\inf} \) lies within it.

\begin{enumerate}[leftmargin=*]
    \item \textbf{Case \( a > 0 \):} This implies that \( \mathbf{k}_{\inf} = -\frac{c}{4a} \leq 0 < a \) lies outside the admissible region $\mathbf{k}_{\inf}\notin(a,\infty)$ for all $c \geq 0$. Since \( \vartheta_2(\mathbf{k}) \) is continuous and has no critical point in $(a, \infty)$, the function \( \vartheta_2(\mathbf{k}) \) is monotonic. 
    To find its minimum, remember that \( d\coloneqq \|k_1 + k_2^\top\|^2 \geq \gamma, \forall c \ge 0 \). Then, using the identity \( \lVert k_1 \rVert^2 + \lVert k_2 \rVert^2 = 2\lVert k_2 \rVert^2 + \tr[\mathcal{C}] \), we can write $\lVert k_1 + k_2^\top \rVert^2 = \sqrt{c + 4\mathbf{k}^2} + 2\mathbf{k} \geq \gamma$.
    Since $d(c) = \sqrt{c + 4\mathbf{k}^2} + 2\mathbf{k}$ is increasing in $c$, we evaluate it at $c = 0$ to get $\sqrt{c + 4\mathbf{k}^2} + 2\mathbf{k} \ge 4\mathbf{k}$, 
    yielding the following inequality $\mathbf{k} \geq \frac{\gamma}{4}$.
    
    Since $\vartheta_2(\mathbf{k})$ is monotonic in $\mathbf{k}$ 
    we substitute $\mathbf{k}$ with $\frac{\gamma}{4}$ to get the uniform upper bound $\vartheta_2(\frac{\gamma}{4})$ while the uniform lower bound over $\mathbf{k} \in (a, \infty)$ comes from the limit such that $\lim_{\mathbf{k} \to a^+} \vartheta_2(\mathbf{k}) = +\infty$ and $\lim_{\mathbf{k} \to \infty} \vartheta_2(\mathbf{k}) = 1$. Hence, we conclude:
    \begin{equation}\label{eq:proof:thm:gPLI:g1}
        g_1(c,\gamma) \coloneqq \sqrt{ \frac{4c + \gamma^2}{(\gamma - 4a)^2} } \ge \vartheta_2(\mathbf{k}) \geq 1,
    \end{equation}
    $\forall \mathbf{k} \in (a, \infty)$. Note that even when $c = 0$, we still have the gPŁI since $\gamma > 4a$ in this case. 
    
    \item \textbf{Case \( a = 0 \):} Using the same argument for the case $a > 0$, we have $4\mathbf{k}> 0$ and $\gamma > 0$. Since $4\mathbf{k} \ge \gamma$ is independent of $a$, then $\mathbf{k} \ge \frac{\gamma}{4}$ is still valid and hence, we conclude:
    \begin{equation}\label{eq:proof:thm:gPLI:g2}
        g_2(c,\gamma) \coloneqq \sqrt{ \frac{4c + \gamma^2}{\gamma^2} } 
        \ge \vartheta_2(\mathbf{k}) \geq  1,
    \end{equation}
     $\forall \mathbf{k} \in (0, \infty)$ and $\forall c \geq 0$, the gPŁI holds since $\gamma > 0$.

    \item \textbf{Case \( a < 0 \):} In this case, the critical point \( \mathbf{k}_{\inf} = -\frac{c}{4a} \geq 0 > a \) lies within the admissible region \( \mathbf{k}_{\inf} \in (a, \infty) \) for all \( c \geq 0 \). To understand the behavior of \( \vartheta_2 \), we evaluate the inner function \( \vartheta_3 \) at \( \mathbf{k}_{\inf} \), resulting $\vartheta_3(\mathbf{k}_{\inf}) = \frac{c}{4a^2 + c}$.
    This expression is strictly positive for all \( c > 0 \), and equals zero if and only if \( c = 0 \), leading to \( \vartheta_2(\mathbf{k}_{\inf}) = 0 \), which invalidates the gPŁI.
    To circumvent this, we again leverage the condition \( d\coloneqq \|k_1 + k_2^\top\|^2 \geq \gamma > 0 \) to refine the lower bound to \( \vartheta_2(\mathbf{k}) > 0 \), as detailed next.

    We now consider a constraint-based lower bound for \( \mathbf{k} \), denoted by \( \underline{\mathbf{k}} \), derived from the condition \( \|k_1 + k_2^\top\|^2 \geq \gamma \) or $d(c) = \sqrt{c + 4\mathbf{k}^2} + 2\mathbf{k} \geq \gamma$, such that
    \begin{equation*}\label{eq:proof:thm:gPLI:k-und}
        \mathbf{k} \geq \frac{\gamma^2 - c}{4\gamma} \eqqcolon \underline{\mathbf{k}}.
    \end{equation*}
    This shows that \( \mathbf{k} \in [\underline{\mathbf{k}}, \infty) \cap (a, \infty) \). When \( \gamma = 0 \), the inequality \( \sqrt{c + 4\mathbf{k}^2} + 2\mathbf{k} \geq 0 \) implies \( c \geq 0 \), and the constraint set \( \mathcal{K}_{\gamma} \) coincides with the admissible set \( \mathcal{K} \coloneqq \{ (k_1, k_2) : \mathbf{k} > a \} \). Thus, for any \( \gamma > 0 \), the definition of \( \mathcal{K}_{\gamma} \) imposes a stricter requirement, ensuring that \( \vartheta_2(\mathbf{k}) > 0 \). Notably, when \( c = 0 \)—a degenerate case in which the critical point \( \mathbf{k}_{\inf} = 0 \) results in \( \vartheta_2(\mathbf{k}_{\inf}) = 0 \)—the constraint \( d \geq \gamma > 0 \) guarantees \( \underline{\mathbf{k}} > 0 \), thereby maintaining strict positivity of \( \vartheta_2(\mathbf{k}) \) throughout the admissible region.
    
    We now examine whether the critical point lies within the admissible region, that is \( \mathbf{k}_{\inf} \in [\underline{\mathbf{k}}, \infty) \), or not. Since \( \vartheta_2(\mathbf{k}) \) has a unique critical point, then if $\mathbf{k}_{\inf} \in [\underline{\mathbf{k}}, \infty)$, we have $\vartheta_2(\mathbf{k}) \geq \vartheta_2(\mathbf{k}_{\inf}), \forall\, \mathbf{k} \in [\underline{\mathbf{k}}, \infty)$. Otherwise, if \( \mathbf{k}_{\inf} < \underline{\mathbf{k}} \), then \( \vartheta_2(\mathbf{k}) \) is monotonically increasing over the admissible region, and we obtain $\vartheta_2(\mathbf{k}) \geq \vartheta_2(\underline{\mathbf{k}}), \forall\, \mathbf{k} \in [\underline{\mathbf{k}}, \infty)$. In both cases, we conclude
    \[
    \vartheta_2(\mathbf{k}) \geq \vartheta_2(\max\{\mathbf{k}_{\inf}, \underline{\mathbf{k}}\}) > 0.
    \]

    We consider two mutually exclusive subcases:
    \begin{itemize}
        \item \textbf{If \( \mathbf{k}_{\inf} \geq \underline{\mathbf{k}} \):} This condition implies
        \[
            -\frac{c}{4a} \geq \frac{\gamma^2 - c}{4\gamma} \quad \Longrightarrow \quad c \geq \frac{a\gamma^2}{a - \gamma} \eqqcolon \Tilde{c} > 0.
        \]
        Since \( a < 0 \) and \( \gamma > 0 \), \( \Tilde{c} \) is strictly positive. Evaluating \( \vartheta_3(\mathbf{k}) \) at  $\mathbf{k}_{\inf}\coloneqq -\frac{c}{4a}$ yields $\vartheta_3(\mathbf{k}_{\inf}) = \frac{c}{4a^2 + c}$
        such that,
        \begin{equation}\label{eq:proof:thm:gPLI:g3}
            \vartheta_2(\mathbf{k}) \geq \sqrt{\frac{c}{4a^2 + c}} \eqqcolon g_3(c) > 0, \quad \forall c \geq \Tilde{c}.
        \end{equation}
        Therefore, under the condition that \( \mathbf{k}_{\inf} \geq \underline{\mathbf{k}} \), the function \( \vartheta_2(\mathbf{k}) \) is uniformly bounded away from zero.
    
        \item \textbf{If \( \underline{\mathbf{k}} \geq \mathbf{k}_{\inf} \):} This condition occurs precisely when \( c \in [0, \Tilde{c}] \). The lower bound of \( \vartheta_2(\mathbf{k}) \) is achieved at \( \underline{\mathbf{k}} \), and thus we analyze the function \( \vartheta_3(\underline{\mathbf{k}}) \). Substituting \( \underline{\mathbf{k}}(c) \coloneqq \frac{\gamma^2 - c}{4\gamma} \) yields
        \begin{equation*}\label{eq:proof:thm:gPLI:v3-k-und}
            \vartheta_3(\underline{\mathbf{k}}(c)) = \frac{(\gamma^2 + c)^2}{(\gamma^2 - c - 4a\gamma)^2}.
        \end{equation*}
        For $c \ge 0$, the numerator $(\gamma^2 + c)^2$ is strictly positive, and at $c=0$ the denominator $(\gamma^2 - 4a\gamma)^2$ is also strictly positive because $a<0$ implies $-4a\gamma>0$. Hence $\vartheta_3(\underline{\mathbf{k}}(c)) > 0$ for all $c \ge 0$ except where the denominator vanishes, when $c = \gamma^2 - 4a\gamma$.
        
        We now show that the interval endpoint $\tilde{c}$ lies strictly to the left of this vertical asymptote, where $\tilde{c} < \gamma^2 - 4a\gamma$. Indeed, the inequality $\frac{a\gamma^2}{a-\gamma} < \gamma^2 - 4a\gamma$ is equivalent (since $a-\gamma < 0$) to
        $\gamma^3 - 4a\gamma^2 + 4a^2\gamma > 0$. Factoring yields $\gamma(\gamma - 2a)^2 > 0$, which is true for all $\gamma>0$ and $a<0$. Therefore, we have $\tilde{c} < \gamma^2 - 4a\gamma$, such that
        \begin{equation}\label{eq:proof:thm:gPLI:g4}
            \vartheta_2(\mathbf{k}) \geq \sqrt{\frac{(\gamma^2 + c)^2}{(\gamma^2 - c - 4a\gamma)^2}} \eqqcolon g_4(c,\gamma) > 0,
        \end{equation}
        for all $c < \Tilde{c}$. In conclusion, $g_4(c,\gamma)$ is strictly positive for $c \in [0,\tilde{c}]$, and the blow-up at $c = \gamma^2 - 4a\gamma$ occurs strictly to the right of $\tilde{c}$, which does not happen since for $c\geq \tilde{c}$ we have $g_3(c)$. Therefore, under the condition that \( \underline{\mathbf{k}} \geq \mathbf{k}_{\inf} \), the function \( \vartheta_2(\mathbf{k}) \) is also uniformly bounded away from zero.
    \end{itemize}
\end{enumerate}

Combining the results of the two terms, $\vartheta_1(\varepsilon)$ and $\vartheta_2(\mathbf{k})$, analyses, we conclude that for any $a$, $c\geq 0$ and $\gamma > \max(0,4a)$, it guarantees the existence of a positive function $\mu_\gamma(c,\gamma) > 0$.
This confirms that the \gPLI holds for the overparameterized LQR problem under the given assumptions. 

    \subsection*{Proof of Corollary~\ref{cor:mu-gPLI}}
    We consider the key results from Theorem~\ref{thm:gPLI}.
    We skip the analysis of $\vartheta_1(\varepsilon)$ since $\vartheta_1(\varepsilon)\ge\frac{r}{4}$ for arbitrary $a$ \eqref{eq:proof:thm:gPLI:v1}. As for $\vartheta_2(\mathbf{k})$, we will discuss \eqref{eq:proof:thm:gPLI:g1}, \eqref{eq:proof:thm:gPLI:g2}, \eqref{eq:proof:thm:gPLI:g3}, and \eqref{eq:proof:thm:gPLI:g4}. 
    
    For $a \geq 0$, notice that $g_1(c,\gamma)$ in \eqref{eq:proof:thm:gPLI:g1} and $g_2(c,\gamma)$ in \eqref{eq:proof:thm:gPLI:g2} are the upper bound functions and $\vartheta_2(\mathbf{k})$ is lower bounded by $\lim_{\mathbf{k} \to \infty} \vartheta_2(\mathbf{k}) = 1$. 
    
    For $a < 0$, we study the two lower-bounding functions, $g_3(c),\forall c \geq \tilde{c}$ and $g_4(c,\gamma)\in[0,\tilde{c}]$. It is obvious that $g_3(c)$ is monotonically increasing in $c$ with minimum at $g_3(\tilde{c})$ and maximum at $\lim_{c\to\infty}g_3(c) = 1$. However, this is not the case for $g_4(c,\gamma)$.
    At \( c = 0 \), we find that the two functions, $g_4(0,\gamma)$ and $g_1(0,\gamma)$, are equivalent, though the values are not because of the set $a$. To determine whether the minimum of \( \vartheta_3(\underline{\mathbf{k}}(c)) \) in the interval \( [0, \Tilde{c}] \) occurs at an interior point or at a boundary, we examine the critical points of \( \vartheta_3(\underline{\mathbf{k}}(c)) \). Taking the derivative $\nabla_{c} \vartheta_3(\underline{\mathbf{k}}(c))$ and solving $\nabla_{c} \vartheta_3(\underline{\mathbf{k}}(c_{\inf})) = 0$ yields the critical point
    \begin{align*}
        2\gamma^3(\gamma - 2a) &= (2a - \gamma)2\gamma c_{\inf}, \quad \longrightarrow \quad
        c_{\inf} &= -\gamma^2
    \end{align*}
    which lies outside the admissible interval \( c_{\inf} \notin [0, \Tilde{c}] \) since $\gamma > 0$. Thus, the minimum of \( \vartheta_3(\underline{\mathbf{k}}(c)) \) over this interval must occur at one of the endpoints. Therefore:
    \begin{equation*}
        \vartheta_2(\mathbf{k}) \geq \min\left\{ g_4(0,\gamma), g_4(\Tilde{c},\gamma) \right\} > 0, \quad \text{for all } c \in [0, \Tilde{c}].
    \end{equation*}
    Finally, comparing $g_4(0,\gamma)$ and $g_4(\tilde{c},\gamma)\equiv g_3(\tilde{c})$:
    \begin{equation*}
        g_4(0,\gamma) = \frac{\gamma^2}{\gamma^2 - 8a\gamma + 16a^2} < \frac{\gamma^2}{\gamma^2 - 4a\gamma + 4a^2} = g_4(\Tilde{c},\gamma),
        \end{equation*}
    implying:
    \begin{equation*}
        \vartheta_2(\mathbf{k}) \geq \sqrt{\frac{\gamma^2}{(\gamma - 4a)^2}} \eqqcolon g_4(0,\gamma) > 0, \quad \text{for all } c \in [0, \Tilde{c}].
    \end{equation*}
    Combining the two sets $c \in [0, \Tilde{c}]$ and $c\ge \tilde{c}$, we obtain a uniform lower bound for \( \vartheta_2(\mathbf{k}) \) over the admissible domain. Notably, when \( c = \Tilde{c} \), the value \( g_3(\Tilde{c}) \) coincides with \( g_4(\Tilde{c},\gamma) \), and it is always bounded below by \( g_4(0,\gamma) \). 
    
    Hence, we conclude \eqref{eq:vec:lowerbound} and
    this establishes the smallest \emph{uniform} value $\underline{\mu}$ of the function $\mu_\gamma(c,\gamma)$ such that the \gPLI holds for arbitrary $a$ and for $(k_1,k_2)\in\mathcal{K}_\gamma$.

    \subsection*{Proof of Corollary~\ref{cor:convergence}}
     We prove the result by comparing the function $\mu_\gamma(c,\gamma)$ by analyzing some arguments developed in the proof of Theorem~\ref{thm:gPLI}. The function $\mu_\gamma(c,\gamma)$ is affected by $\vartheta_1(\varepsilon)$, bounded below by $\frac{r}{4}$ (independent of $c$), and $\vartheta_2(\mathbf{k})$ using $g_i(c),\forall i = 1,\dots,4$. These functions $g_i(c)>0,\forall i$ given in \eqref{eq:proof:thm:gPLI:g1}, \eqref{eq:proof:thm:gPLI:g2}, \eqref{eq:proof:thm:gPLI:g3}, and \eqref{eq:proof:thm:gPLI:g4} are monotonically increasing in $c$ for a fixed $\gamma$.
     This concludes that for the two initializations  \( (\bar{k}_1, \bar{k}_2), (\hat{k}_1, \hat{k}_2) \in \mathcal{K}_{\gamma} \) with the same cost \( \mathcal{L}(\bar{k}_1, \bar{k}_2) = \mathcal{L}(\hat{k}_1, \hat{k}_2) \) and different imbalance \( \bar{c} > \hat{c} \), we have \( \mu_{\gamma}(\bar{c},\gamma) > \mu_{\gamma}(\hat{c},\gamma) \), resulting:
     \begin{equation*}
        e^{-\mu_{\gamma}(\bar{c},\gamma) t} \big(\mathcal{L}(\bar{k}_1, \bar{k}_2) - \underline{\mathcal{L}}\big) < e^{-\mu_{\gamma}(\hat{c},\gamma) t} \big(\mathcal{L}(\hat{k}_1, \hat{k}_2) - \underline{\mathcal{L}}\big).
    \end{equation*}
    This completes the proof.

\end{document}